\documentclass[preprint,sort&compress,final,12pt]{elsarticle}
\usepackage{graphicx}
\usepackage{epstopdf}
\usepackage{amsmath}
\usepackage{amssymb}
\usepackage{latexsym}
\usepackage{natbib}
\usepackage{color}
\usepackage{hyperref}
\usepackage{geometry}
\usepackage{bm}
\usepackage{indentfirst}

\usepackage{cases}
\usepackage{datetime}
\newtheorem{theorem}{Theorem}[section]

\newtheorem{assumption}{Assumption}[section]

\newtheorem{lemma}[theorem]{Lemma}
\newtheorem{proposition}[theorem]{Proposition}
\newtheorem{remark}[theorem]{Remark}

\numberwithin{equation}{section}

\def\Proof{\noindent{\bf Proof.}~}

\def\Re{\mathrm{Re}}
\def\Im{\mathrm{Im}}

\journal{\empty}
\date{\empty}
\begin{document}
	
\begin{frontmatter}
		
\title{A sufficient condition on successful invasion by the predator}

\author{Hongliang Li\footnote{Corresponding author.}}
\ead{lihl@mail.bnu.edu.cn}

\author{Min Zhao}
\ead{minzhao@mail.bnu.edu.cn}

\author{Rong Yuan}
\ead{ryuan@bnu.edu.cn}

\address{Laboratory of Mathematics and Complex Systems (Ministry of Education), School of Mathematical Sciences, Beijing Normal University, Beijing 100875, People's Republic of China}
		
\begin{abstract}
In this paper, we provide a sufficient condition on successful invasion by the predator. Specially, we obtain the persistence of traveling wave solutions of predator-prey system, in which the predator can survive without the predation of the prey. This proof heavily depends on comparison principle of scalar monostable equation, the rescaling method and phase-plane analysis.
\end{abstract}
		

\end{frontmatter}

\section{Introdcution}
In the natural world, we are more concerned with the ideal scenario in which the predator and the prey may cohabit.
As a consequence, many researchers investigate the process of population ecology propagation, which may be well described by traveling wave solutions connecting the prey-present state and the coexistence state, also named invasion waves.
At present, the method of dealing with the existence of traveling wave solutions, in general, includes the shoot argument \cite{dunbar1983travelling}, Schauder's fixed-point theorem\cite{ma2001traveling} and Conley index \cite{gardner1984existence}. It is generally accepted that the asymptotic behavior of traveling wave solutions at positive infinity often depends on an appropriate and technicality Lyapunov function, nevertheless, construction of the Lyapunov function is not easy, particularly for the non-monotonic functional response. Just as said in \cite{zhang2016minimal}, it is sufficient to study the persistence of traveling wave solutions if we only want to know whether the invasion is successful and what the invasion speed is.

Therefore, the purpose of this work is to investigate the persistence of traveling wave solutions to general predator-prey system with diffusion
\begin{equation}\label{eq:g1}
\left\{\begin{split}
&u_t= u_{xx}+uF(u,v),\\[0.2cm]
&v_t=d v_{xx}+vG(u,v),
\end{split}
\right.
\end{equation}
where $u$ and $v$ stand for the population densities of the prey and the predator, respectively. $d>0$ denotes the ratio of the diffusion of the predator to that of the prey, $F(u,v)$ and $G(u,v)$ describe the interaction between the predator and the prey, and the assumptions for these are given below.

\begin{assumption}\label{a:1}
$F(u,v)$ and $G(u,v)$ are $C^2$ functions and satisfy
\begin{description}
  \item(a) There is $\mu$ such that $G(0,v)(v-\mu)<0$  for  $v\in[0,\mu)\cup(\mu,+\infty)$.
  \item(b) There is $v_0$ such that $G(u,v_0)<0$ for $u\in[0,1]$, and
  \begin{equation*}
    G(1,0)\geq G(u,v)\geq G(0,v) \text{ for } (u,v)\in[0,1]\times[0,v_0].
  \end{equation*}
  \item(c) $F(1,0)=0$, $F(u,0)>0$ for $u\in[0,1)$, and $F(1, v)<0$ for $v\in(0,v_0]$.
\end{description}
\end{assumption}

In view of Assumptions \ref{a:1}, $u$-equation and $v$-equation are reduced to a scalar monostable reaction-diffusion equation in the absence of the predator and the prey, respectively. And as for the corresponding kinetic system, there are always three boundary equilibria $(1,0)$, $(0,0)$ and $(0,\mu)$. Some classical examples satisfying the above assumptions are as follows
\begin{description}
  \item (1) Lotka-Volterra predator-prey system: for given $a,r,b>0$
  \begin{equation*}
    F(u,v)=1-u-av \text{ and } G(u,v)=r(1-v+bu).
  \end{equation*}
  \item (2) Modified Leslie-Gower predator-prey system: for given $a, e_1, r, e_2>0$
    \begin{equation*}
    F(u,v)=1-u-\displaystyle\frac{av}{u+e_1} \text{ and } G(u,v)=r\left(1-\displaystyle\frac{v}{u+e_2}\right).
  \end{equation*}
\end{description}
A positive solution is called a traveling wave solution if solution of \eqref{eq:g1} has the form
\begin{equation*}
  (u,v)(x,t)=(\phi,\psi)(z), \ z=x+ct,
\end{equation*}
where $c$ is the wave speed. Automatically, $(\phi,\psi)(z)$ satisfies
\begin{equation}\label{eq:tr}
\left\{\begin{split}
&\phi''(z)-c\phi'(z)+\phi(z)F(\phi,\psi)(z)=0,\\[0.2cm]
&d\psi''(z)-c\psi'(z)+\psi(z)G(\phi,\psi)(z)=0.
\end{split}
\right.
\end{equation}
Traveling wave solution is persistent if there exist two positive constants $M_1$ and $M_2$ such that
\begin{equation*}
\left\{\begin{array}{ll}
(\phi, \psi)(-\infty)=(1, 0), \\[0.2cm]
M_{1}<\liminf _{z \rightarrow +\infty} \phi(z)\leq\limsup _{z \rightarrow +\infty} \phi(z)<M_{2}, \\[0.2cm]
M_{1}<\liminf _{z \rightarrow +\infty} \psi(z)\leq\limsup _{z \rightarrow +\infty} \psi(z)<M_{2}.
\end{array}
\right.
\end{equation*}

In this paper, we always assume that $c\geq c_*:=2\sqrt{dG(1,0)}$ and define
\begin{equation*}
   \lambda_1:=\frac{1}{2d}\left(c-\sqrt{c^2-4dG(1, 0)}\right)>0.
\end{equation*}
Now, we state the main result.
\begin{theorem}\label{th:ma}
For $c\geq c_*$, if $F(0,\mu)>0$, then traveling wave solutions of system \eqref{eq:g1} are persistent.
\end{theorem}

The proof of Theorem \ref{th:ma} is given in Section 2. Firstly, we obtain the existence of a positive solution of system \eqref{eq:tr} from Theorem 2.1 in \cite{ai2017traveling}. Then based on the spreading theory on the scalar monostable equation and comparison principle, we obtain that component $\psi(z)$ of traveling wave solutions is persistent. Meanwhile, inspired by \cite{yang2020wave}, we show that component $\phi(z)$ of traveling wave solutions is persistent by using the rescaling method and phase-plane analysis.

\section{The proof of Theorem 1.1}
Let us begin this section with the existence result, which may be deduced immediately from Theorem 2.1 in \cite{ai2017traveling}.
\begin{lemma}\label{th:1}
For $c\geq c_*$, system \eqref{eq:tr} admits a positive solution $(\phi,\psi)(z)$ satisfying
\begin{equation*}
  (\phi,\phi',\psi,\psi')(-\infty)=(1,0,0,0) \text{ and } 0<\phi(z)<1, \ 0<\psi(z)<v_0 \text{ over } \mathbb R.
\end{equation*}
Moreover, $\phi'(z)$ and $\psi'(z)$ are bounded over $\mathbb R$.
\end{lemma}

To proceed, we require a priori estimate on the inhomogeneous linear equation
\begin{equation}\label{eq:nl}
  w''(z)+\alpha w'(z)+f(z)w(z)=h(z),
\end{equation}
where $\alpha$ is a positive constant and $f$, $h$ are continuous functions on $[a, b]$.
\begin{proposition}\label{pr:1}
(Lemma 3.3 of \cite{fu2014existence})
Assume that $w(z)\in C([a,b])\cap C^2((a,b))$ satisfies \eqref{eq:nl} in $(a,b)$ and
\begin{equation*}
  -\gamma_1\leq f\leq0 \text{ and } |h|\leq \gamma_2 \text{ on } [a,b],
\end{equation*}
for some positive constants $\gamma_i$, $i=1,2$. If $||w||_{C([a,b])}\leq \gamma_3$ for positive constant $\gamma_3$, then there exists a positive constant $\gamma_4(\alpha, \gamma_1, \gamma_2, \gamma_3, a, b)$  such that
\begin{equation*}
  ||w'||_{C([a,b])}\leq \gamma_4(\alpha, \gamma_1, \gamma_2, \gamma_3, a, b).
\end{equation*}
\end{proposition}

The following two lemmas describe the persistence of travelling wave solutions.
\begin{lemma}\label{le:v}
$\liminf_{z\rightarrow+\infty} \psi(z)\geq \mu$.
\end{lemma}
\Proof From Lemma \ref{th:1}, we have $\psi(z)>0$ over $\mathbb R$ and $\psi(-\infty)=0$ for $c\geq c_*$, then there is a constant $\zeta\in(0,\mu)$ such that $\psi(z_1)=\zeta$ for $z_1\in\mathbb R$. Since $(\phi,\psi)(z+a)$ is also traveling wave solution for any $a\in\mathbb R$, we consider $z_1=0$ without losing generality. Note that $\psi(z)$ is uniformly continuous over $\mathbb R$ owing to the boundedness of $\psi'(z)$, there is a $\epsilon>0$ such that $\psi(z)>\zeta/2$ for $z\in(-\epsilon,\epsilon)$.
Now, we consider initial value problem
\begin{equation}\label{eq:ms}
\left\{\begin{split}
&w_t= d\Delta w+wg(w),\\[0.2cm]
&w(0,x)=\varphi(x),
\end{split}
\right.
\end{equation}
where $g(w):=G(0,w)$ and uniformly continuous function $\varphi(x)$ satisfies the conditions:
\begin{description}
  \item(1) $\varphi(x)=\zeta/2$ for $x\in[-\epsilon/2,\epsilon/2]$ and $\varphi(x)=0$ for $\mathbb R\setminus (-\epsilon,\epsilon)$,
  \item(2) $\varphi(x)$ is increasing for $x\in[-\epsilon,-\epsilon/2]$ and decreasing for  $x\in[\epsilon/2,\epsilon]$.
\end{description}
From Assumption \ref{a:1}-\emph{(a)}, we know that
\begin{equation*}
   g(\mu)=0 \text{ and } g(v)>0  \text{ for }  v\in[0,\mu).
\end{equation*}
According to Proposition 3.1 in \cite{ducrot2019spreading}, there exists a $c^*(d,g(v))\geq2\sqrt{dg(0)}$ such that
\begin{equation*}
 \lim\limits_{t\rightarrow+\infty}\inf\limits_{|x|<c^*(d,g(v))t}w(x,t)= \mu.
\end{equation*}
Thanks to $\phi(z)>0$ and Assumption \ref{a:1}-\emph{(b)}, then $v(x,t):=\psi(z)$ satisfies
\begin{equation*}
\left\{\begin{split}
  &v_t\geq d\Delta v+vg(v),\\[0.2cm]
  &v(x,0)=\psi(x).
\end{split}
\right.
\end{equation*}
Obviously, $\psi(x)\geq\varphi(x)$ over $\mathbb R$. Comparison principle suggests that for $c\geq c_*$
\begin{equation*}
  \liminf\limits_{z\rightarrow+\infty}\psi(z)=\liminf\limits_{z\rightarrow+\infty}v(0,z/c)\geq \liminf\limits_{z\rightarrow+\infty}w(0,z/c)\geq \lim\limits_{t\rightarrow+\infty}\inf\limits_{|x|<c^*(d,g(v))t}w(x,t)
  = \mu.
\end{equation*}
Hence, we complete the proof.

\begin{lemma}\label{le:2}
$\liminf_{z\rightarrow+\infty} \phi(z)>0$ if $F(0,\mu)>0$.
\end{lemma}
\Proof For contradiction, we assume that there exists a sequence $\left\{z_n\right\}_{n\in\mathbb N}$ satisfying $z_n\rightarrow+\infty$ as $n\rightarrow+\infty$ such that $\phi(z_n)\rightarrow0$ as $n\rightarrow+\infty$. Let $\left(\phi_n,\psi_n\right)(z):=(\phi, \psi)(z_n+z)$, then $\phi''_n(z)$ and $\psi''_n(z)$ are bounded over $\mathbb R$ from Lemma \ref{th:1}. Differentiating system \eqref{eq:tr}, we can express $\phi'''_n(z)$ and $\psi'''_n(z)$ in terms of $\phi_n(z), \phi'_n(z), \phi''_n(z), \psi_n(z), \psi'_n(z)$ and $\psi''_n(z)$, hence $\phi'''_n(z)$ and $\psi'''_n(z)$ are bounded over $\mathbb R$. Therefore, by Arzel\`{a}-Ascoli theorem, up to extracting a subsequence, there are some functions $\phi_{\infty}(z)$ and $\psi_{\infty}(z)$ such that
$\phi_n(z)\rightarrow \phi_{\infty}(z)$ and $\psi_n(z)\rightarrow \psi_{\infty}(z)$ in $C_{loc}^2(\mathbb R)$  as $n\rightarrow+\infty$.
Moreover, $\phi_{\infty}(z)$ satisfies $0\leq\phi_{\infty}(z)\leq1$ over $\mathbb R$ and
\begin{align*}
  \phi''_\infty(z)-c\phi'_\infty(z)+\phi_\infty(z)F(\phi_\infty,\psi_\infty)(z)=0.
\end{align*}
Note that $\phi_{\infty}(0)=0$ due to $\phi(z_n)\rightarrow0$ as $n\rightarrow+\infty$, we have $\phi_{\infty}(z)\equiv0$ over $\mathbb R$. Hence, $\psi_{\infty}(z)$ should satisfy $0\leq\psi_{\infty}(z)\leq v_0$ over $\mathbb R$ and
\begin{align*}
  d\psi''_\infty(z)-c\psi'_\infty(z)+\psi_\infty(z)G(0,\psi_\infty(z))=0.
\end{align*}
It follows from Lemma \ref{le:v} that $\psi_{\infty}(z)\geq\mu$ for $z\geq0$, one can readily check that $0<\psi_{\infty}(z)\leq v_0$ over $\mathbb R$.
We will show that $(\psi_{\infty},\psi'_\infty)(0)=(\mu,0)$ by phase-plane analysis. In order to do it,
we set $\chi_\infty(z):=\psi'_\infty(z)$ and rewrite $\psi_{\infty}$-equation as a system of first order ODEs in $\mathbb{R}^2$,
\begin{equation}\label{eq:vv}
\left\{\begin{split}
  &\psi'_\infty(z)=\chi_\infty(z),\\[0.2cm]
  &d\chi_\infty'(z)=c\chi_\infty(z)-\psi_{\infty}(z)G(0,\psi_{\infty}(z)).
\end{split}
\right.
\end{equation}
Moreover, some sets are given by
\begin{align*}
  A_1&:=\left\{(\psi_{\infty},\chi_\infty)(z): \psi_{\infty}(z)>\mu, \chi_\infty(z)>0\right\},\\[0.2cm]
  A_2&:=\left\{(\psi_{\infty},\chi_\infty)(z): \psi_{\infty}(z)=\mu, \chi_\infty(z)>0\right\},\\[0.2cm]
  A_3&:=\left\{(\psi_{\infty},\chi_\infty)(z): \psi_{\infty}(z)>\mu, \chi_\infty(z)=0\right\},\\[0.2cm]
  A_4&:=\left\{(\psi_{\infty},\chi_\infty)(z): \psi_{\infty}(z)=\mu, \chi_\infty(z)<0\right\},\\[0.2cm]
  A_5&:=\left\{(\psi_{\infty},\chi_\infty)(z): \psi_{\infty}(z)>\mu, \chi_\infty(z)<0\right\},\\[0.2cm]
  A_6&:=\left\{(\psi_{\infty},\chi_\infty)(z): \psi_{\infty}(z)=\mu, \chi_\infty(z)=0\right\}.
\end{align*}
If there exists a $z_0\in\mathbb R$ such that $(\psi_{\infty},\chi_\infty)(z_0)\in A_1$, one can see that $\psi_{\infty}(+\infty)=+\infty$, which contradicts to the boundedness of $\psi_{\infty}(z)$. Hence we have $(\psi_{\infty},\chi_\infty)(z)\notin A_1$ over $\mathbb R$.
Meanwhile, if there exists a $z_0\in\mathbb R$ such that $(\psi_{\infty},\chi_\infty)(z_0)\in A_2\cup A_3$, obviously, we have $(\psi_{\infty},\chi_\infty)(z^+_0)\in A_1$. Thus, $(\psi_{\infty},\chi_\infty)(z)\notin A_1\cup A_2\cup A_3$ over $\mathbb R$. Undoubtedly, $(\psi_{\infty},\chi_\infty)(0)\notin A_1\cup A_2\cup A_3$. Let us continuous to show that $(\psi_{\infty},\chi_\infty)(0)\notin A_4\cup A_5$.
If $(\psi_{\infty},\chi_\infty)(0)\in A_4$, clearly, $\psi_{\infty}(0^+)<\mu$, which contradicts to $\psi_{\infty}(z)\geq\mu$ for $z\geq0$. On the other hand, if $(\psi_{\infty},\chi_\infty)(0)\in A_5$, then $\chi_\infty(z)<0$ for $z\leq0$ owing to $(\psi_{\infty},\chi_\infty)(z)\notin A_3$ over $\mathbb R$, which suggests that for $z\leq0$
\begin{align*}
  (\psi_{\infty},\chi_\infty)(z)\in\left\{(\psi_{\infty},\chi_\infty)(z): v_0\geq \psi_{\infty}(z)>\mu, \chi_\infty(z)<0\right\}.
\end{align*}
Thus, monotone bounded principle explains that $\psi_{\infty}(-\infty)$ exists and $\psi_{\infty}(-\infty)>\psi_{\infty}(0)\geq\mu$.
In the meantime, we claim that $\vartheta_-\leq\chi_\infty(z)\leq\vartheta_+$ over $\mathbb R$, where
\begin{equation*}
  \vartheta_-:=\min_{\eta\in[0,v_0]}\left\{\eta G(0,\eta)/c\right\}<0<\vartheta_+:=\max_{\eta\in[0,v_0]}\left\{\eta G(0,\eta)/c\right\}
\end{equation*}
from Assumptions \ref{a:1}-\emph{(a)}. We only show the inequality $\chi_\infty(z)\leq\vartheta_+$ here since the inequality $\chi_\infty(z)\geq\vartheta_-$ then can be obtain in same way. For contradiction, we assume that there exists a $z_0$ such that $\chi_\infty(z_0)>\vartheta_+$, then one can easily see that $\chi_\infty(z)>\vartheta_+$ for $z>z_0$ from \eqref{eq:vv}, which contradicts to the boundedness of $\psi_\infty(z)$. This claim gives the boundedness of $\chi'_\infty(z)$, and
differentiating $\chi_\infty$-equation further yields the boundedness of $\chi''_\infty(z)$,
hence $(\chi_\infty,\chi'_\infty)(-\infty)=0$ by Barbalat Lemma, which suggests that $G(0,\psi_{\infty}(-\infty))=0$. However, this is impossible since $G(0,v)(v-\mu)<0$ for $v>\mu$ in Assumption \ref{a:1}-\emph{(a)}. Therefore, there must be $(\psi_{\infty},\chi_\infty)(0)\in A_6$. Note that $A_6$ is an invariant set of system \eqref{eq:vv}, then $\psi_{\infty}(z)\equiv\mu$ over $\mathbb R$. Now, we define
\begin{equation}\label{eq:s5}
  \omega_n(z)=\frac{\phi_n(z)}{\phi(z_n)}=\frac{\phi(z+z_n)}{\phi(z_n)}=\exp\left\{\int_{z_n}^{z_n+z}\frac{\phi'(z)}{\phi(z)} dz\right\}.
\end{equation}
Clearly, $\omega_n(z)>0$ over $\mathbb R$. It follows from $\phi_n$-equation that
\begin{equation}\label{eq:v}
  \omega_n''(z)-c \omega_n'(z)+\omega_n(z)F(\phi_n,\psi_n)(z)=0.
\end{equation}
We claim that function $\varpi(z)=\phi'(z)/\phi(z)$ is bounded over $\mathbb R$. In fact, since $\phi(z)$ and $\psi(z)$ are bounded, there exists a $m>0$ such that $F(\phi,\psi)(z)\geq -m$ over $\mathbb R$. Automatically,
\begin{align*}
  \varpi'(z)&=c\varpi(z)-\varpi^2(z)-F(\phi,\psi)(z)\\[0.2cm]
  &\leq c\varpi(z)-\varpi^2(z)+m
  =-(\varpi(z)-\pi_+)(\varpi(z)-\pi_-),
\end{align*}
where
\begin{equation*}
  \pi_+=\frac{1}{2}\left(c+\sqrt{c^2+4m}\right) \text{ and } \pi_-=\frac{1}{2}\left(c-\sqrt{c^2+4m}\right).
\end{equation*}
Owing to $\varpi(-\infty)=0$, there exists a $z_0\ll-1$ such that $\varpi(z)<\pi_+$ for $z\in(-\infty, z_0]$. Since function $\varpi_+(z):=\pi_+$ satisfies $\varpi_+'(z)=-(\varpi_+(z)-\pi_+)(\varpi(z)_+-\pi_-)$,
we have $\varpi(z)\leq\varpi_+(z)$ for $z\in(z_0,+\infty]$ by comparison principle. Hence $\varpi(z)\leq\pi_+$ over $\mathbb R$.
Let us continuous to show $\varpi(z)>-\pi_+$ over $\mathbb R$. For contradiction, we assume that $\varpi(z_1)\leq-\pi_+$ for $z_1\in\mathbb R$. Let $\varpi_-(z)$ be solution of the following equation
\begin{equation}\label{eq:g2}
\left\{\begin{split}
  &\varpi_-'(z)=c\varpi_-(z)-\varpi_-^2(z)+m,\\[0.2cm]
  &\varpi_-(z_1)=\varpi(z_1).
\end{split}
\right.
\end{equation}
By solving \eqref{eq:g2}, we get
\begin{align*}
\varpi_-(z)=\frac{\pi_{+}-\pi_{-} e^{-\sqrt{c^{2}+4 m}\left(z-z_{2}\right)}}{1-e^{-\sqrt{c^{2}+4 m}\left(z-z_{2}\right)}},
\end{align*}
where
\begin{equation*}
  z_{2}=z_{1}+\frac{1}{\sqrt{c^{2}+4 m}} \ln \left(\frac{\pi_{+}-\varpi_-\left(z_{1}\right)}{\pi_{-}-\varpi^-\left(z_{1}\right)}\right)>z_{1}.
\end{equation*}
Hence $\varpi_-(z)\rightarrow-\infty$ as $z\rightarrow z_2^-$. On the other hand,
comparison principle yields that $\varpi_-(z)\geq\varpi(z)$ for $z>z_1$, which implies that $\varpi(z)\rightarrow-\infty$ as $z\rightarrow z_3^-$ with $z_3\in(z_1,z_2]$. However, $\varpi(z)$ is defined in all $z\in\mathbb R$, hence the claim is valid. According to the above claim, $\omega_n(z)$ is locally uniformly bounded. Let $\upsilon_n(z):=\omega_n(-z)$, then for any positive constant $\beta$
\begin{equation*}
  \upsilon_n''(z)+c \upsilon_n'(z)-\beta\upsilon_n(z)=-\upsilon_n(z)\left(F(\phi_n,\psi_n)(-z)+\beta\right).
\end{equation*}
Hence $\upsilon_n'(z)$ is locally uniformly bounded by using Proposition \ref{pr:1}, in other words, $\omega_n'(z)$ is locally uniformly bounded. $\omega_n$-equation further gives that $\omega''_n(z)$ is locally uniformly bounded, so is its derivative by differentiating \eqref{eq:v}. Therefore, with the aid of Arzel\`{a}-Ascoli theorem, up to extracting a subsequence, there is a function $\omega_\infty(z)$ such that $\omega_n(z)\rightarrow\omega_\infty(z)$ in $C_{loc}^2(\mathbb R)$ as $n\rightarrow+\infty$. Hence we obtain a linear equation
\begin{equation}\label{eq:u}
  \omega_\infty''(z)-c\omega_\infty'(z)+F(0,\mu)\omega_\infty(z)=0
\end{equation}
due to $(\phi_\infty, \psi_\infty)(z)\equiv(0,\mu)$ over $\mathbb R$. Note that $\omega_n(z)>0$ over $\mathbb R$, we have $\omega_\infty(z)\geq0$ over $\mathbb R$, meanwhile, $\omega_\infty(z)>0$ over $\mathbb R$ thanks to $\omega_\infty(0)=1$.
Firstly, let us consider $c_0:=2\sqrt{F(0,\mu)}>c$. At the moment, we have
\begin{equation*}
  \omega_\infty(z)=e^{z\Re(\eta^+)}\left[cos\left(z\Im(\eta^+)\right)+\kappa sin\left(z\Im(\eta^+)\right)\right] \text{ with } \kappa\in\mathbb R,
\end{equation*}
where
\begin{equation*}
  \eta^\pm=\frac{1}{2}\left(c\pm\sqrt{c^2-4F(0,\mu)}\right).
\end{equation*}
Hence, $\omega_\infty(z_0)<0$ for some $z_0\in\mathbb R$, and we get a contradiction with the positivity of $\omega_\infty(z)$.
On the other hand, we consider $c\geq c_0$. At the moment, there are constants $\kappa_1$ and $\kappa_2$ such that
\begin{align*}
  \omega_\infty(z)=
  \left\{
  \begin{array}{ll}
    \kappa_1 e^{\eta^-z}+(1-\kappa_1)e^{\eta^+z}, &\text{ if } c>c_0, \\[0.2cm]
    e^{\eta^+z}(1+\kappa_2z), &\text{ if } c=c_0.
  \end{array}
  \right.
\end{align*}
Using the positivity of $\omega_\infty(z)$ again, we obtain $\kappa_1\in[0,1]$ and $\kappa_2=0$. And
\begin{align*}
  \omega'_\infty(z)=
  \left\{
  \begin{array}{ll}
    \kappa_1\eta^-e^{\eta^-z}+(1-\kappa_1)\eta^+e^{\eta^+z}, &\text{ if } c>c_0, \\[0.2cm]
    \eta^+e^{\eta^+z}, &\text{ if } c=c_0.
  \end{array}
  \right.
\end{align*}
Therefore, we have $\omega'_\infty(z)>0$ over $\mathbb R$, which implies that $\omega'_n(z)>0$ in any bounded closed set for large enough $n$. However, since $\phi(-\infty)=1$ and $0<\phi(z)<1$ over $\mathbb R$, then there is a $z^*\in\mathbb R$ such that $\phi'(z)\leq0$ for $z\leq z^*$, it follows from \eqref{eq:s5} that $\omega'_n(z)\leq0$ for $z\in[z^*-1, z^*]$. Thus, we get a contradiction, which finish proof of this lemma.

\begin{remark}
If $F(0,\mu)=0$, then $\omega_\infty(z)=\kappa+(1-\kappa)e^{cz}$ for $\kappa\in[0,1]$ owing $\omega_\infty(0)=1$ and $\omega_\infty(z)>0$ over $\mathbb R$. Clearly, $\omega'_\infty(z)>0$ or $\omega_\infty(z)\equiv1$ over $\mathbb R$. We are unable to judge whether traveling wave solutions are persistent for $c\geq c_*$ since the latter cannot be excluded.
\end{remark}

Finally, we finish the proof of Theorem \ref{th:ma}.
\begin{lemma}
$\limsup _{z \rightarrow +\infty} \phi(z)<1$ and $\limsup _{z \rightarrow +\infty} \psi(z)<v_0$.
\end{lemma}
\Proof For contradiction, we assume that there exists a sequence $\left\{z_n\right\}_{n\in\mathbb N}$ satisfying $z_n\rightarrow+\infty$ as $n\rightarrow+\infty$ such that $\phi(z_n)\rightarrow1$ as $n\rightarrow+\infty$. Let $\left(\phi_n,\psi_n\right)(z):=(\phi, \psi)(z_n+z)$, similar to proof in Lemma \ref{le:2}, there are some functions $\phi_{\infty}(z)$ and $\psi_{\infty}(z)$ satisfy
\begin{equation}\label{eq:f}
\begin{split}
  &\phi''_\infty(z)-c\phi'_\infty(z)+\phi_\infty(z)F(\phi_\infty,\psi_\infty)(z)=0,\\[0.2cm]
  &d\psi''_\infty(z)-c\psi'_\infty(z)+\psi_\infty(z)G(\phi_\infty,\psi_\infty)(z)=0.
\end{split}
\end{equation}
One can readily verify that $0<\phi_\infty\leq1$ and $0<\psi_\infty\leq v_0$ over $\mathbb R$. Note that $\phi_{\infty}(0)=1$ due to $\phi(z_n)\rightarrow1$ as $n\rightarrow+\infty$, we have $\phi'_{\infty}(0)=0$ and $\phi''_{\infty}(0)\leq 0$. Thus, from Assumption \ref{a:1}-\emph{(c)},
\begin{align*}
  0\geq \phi''_\infty(0)=-F(1,\psi_\infty(0))>0,
\end{align*}
hence we conclude $\limsup _{z \rightarrow +\infty} \phi(z)<1$. A similar discussion yields that
\begin{align*}
  0\geq d\psi''_\infty(0)=-v_0G(\phi_\infty(0),v_0).
\end{align*}
From Assumption \ref{a:1}-\emph{(b)}, we have $\phi_\infty(0)>1$, which leads to a contradiction. Hence we also conclude $\limsup _{z \rightarrow +\infty} \psi(z)<v_0$. Therefore, we complete the proof.

\begin{remark}
Our result can also be applied to Holling-Tanner system in \cite{zhao2022traveling}
\begin{equation}\label{eq:ht1}
\left\{\begin{split}
&u_t= u_{xx}+u\left(1-u-\frac{\alpha u^{m-1}v}{1+\beta_1 u^m+\beta_2 v^m}\right),\\[0.2cm]
&v_t=d v_{xx}+rv\left(1-\frac{v}{u}\right),
\end{split}
\right.
\end{equation}
where $\alpha, d,r>0$, $\beta_1,\beta_2\geq0$ and $m\geq1$. We follow the idea of \cite{ai2017traveling} to consider
\begin{equation}\label{eq:ht2}
\left\{\begin{split}
&u_t= u_{xx}+u\left(1-u-\frac{\alpha u^{m-1}v}{1+\beta_1 u^m+\beta_2 v^m}\right),\\[0.2cm]
&v_t=d v_{xx}+rv\left(1-\frac{v}{\sigma_\varepsilon(u)}\right),
\end{split}
\right.
\end{equation}
where
\begin{align*}
  \sigma_\varepsilon(u)=\left\{
  \begin{array}{ll}
    u, & u\geq\varepsilon, \\[0.2cm]
    u+\varepsilon e^\frac{1}{u-\varepsilon}, & 0\leq u<\varepsilon.
  \end{array}
  \right.
\end{align*}
Clearly, $\mu=\varepsilon (1/e)^{1/\varepsilon}$ and
\begin{equation*}
  F(0,\mu)=\left\{
  \begin{array}{ll}
  1, & m>1,\\[0.2cm]
  1-\frac{\alpha \varepsilon (1/e)^{1/\varepsilon}}{1+\beta_2 \varepsilon (1/e)^{1/\varepsilon}}, & m=1.
  \end{array}
  \right.
\end{equation*}
Thus, $F(0,\mu)>0$ as long as $\varepsilon$  is small enough. In view of Theorem \ref{th:ma}, we conclude that traveling wave solution $(\phi_\varepsilon,\psi_\varepsilon)(z)$ of system  \eqref{eq:ht2} is persistent for $c\geq2\sqrt{dr}$,  then there exists $\delta>0$ such that $\phi_\varepsilon(z)>\delta$ over $\mathbb R$.
Consequently, $\sigma_\varepsilon(\phi_\varepsilon)=\phi_\varepsilon$ for small enough $\varepsilon$, and $(\phi, \psi)(z):=(\phi_\varepsilon,\psi_\varepsilon)(z)$ is
a traveling wave solution of system \eqref{eq:ht1}.
\end{remark}
\section*{Acknowledgments}
This work is supported by the National Natural Science Foundation of China (No. 12171039 and 12271044).
\section*{Statements and Declarations}
Authors have no conflict of interest to declare.
\bibliographystyle{abbrv}
\bibliography{CKD}

\begin{thebibliography}{1}

\bibitem{ai2017traveling}
S.~Ai, Y.~Du, and R.~Peng.
\newblock Traveling waves for a generalized {Holling-Tanner} predator-prey
  model.
\newblock {\em Journal of Differential Equations}, 263(11):7782--7814, 2017.

\bibitem{ducrot2019spreading}
A.~Ducrot, T.~Giletti, and H.~Matano.
\newblock Spreading speeds for multidimensional reaction-diffusion systems of
  the prey-predator type.
\newblock {\em Calculus of Variations and Partial Differential Equations},
  58:1--34, 2019.

\bibitem{dunbar1983travelling}
S.~R. Dunbar.
\newblock Travelling wave solutions of diffusive {Lotka-Volterra} equations.
\newblock {\em Journal of Mathematical Biology}, 17:11--32, 1983.

\bibitem{fu2014existence}
S.~Fu.
\newblock The existence of traveling wave fronts for a reaction-diffusion
  system modelling the acidic nitrate-ferroin reaction.
\newblock {\em Quarterly of Applied Mathematics}, 72(4):649--664, 2014.

\bibitem{gardner1984existence}
R.~A. Gardner.
\newblock Existence of travelling wave solutions of predator-prey systems via
  the connection index.
\newblock {\em SIAM Journal on Applied Mathematics}, 44(1):56--79, 1984.

\bibitem{ma2001traveling}
S.~Ma.
\newblock Traveling wavefronts for delayed reaction-diffusion systems via a
  fixed point theorem.
\newblock {\em Journal of Differential Equations}, 171(2):294--314, 2001.

\bibitem{yang2020wave}
F.~Yang, W.~Li, and J.~Wang.
\newblock Wave propagation for a class of non-local dispersal non-cooperative
  systems.
\newblock {\em Proceedings of the Royal Society of Edinburgh Section A:
  Mathematics}, 150(4):1965--1997, 2020.

\bibitem{zhang2016minimal}
T.~Zhang, W.~Wang, and K.~Wang.
\newblock Minimal wave speed for a class of non-cooperative diffusion-reaction
  system.
\newblock {\em Journal of Differential Equations}, 260(3):2763--2791, 2016.

\bibitem{zhao2022traveling}
X.~Zhao and H.~Wang.
\newblock Traveling waves for a generalized {Beddington-DeAngelis}
  predator-prey model.
\newblock {\em Communications in Nonlinear Science and Numerical Simulation},
  111:106478, 2022.

\end{thebibliography}

\end{document}